# Why do we change whatever amount we found in the first envelope: the Wikipedia "two envelopes problem" commented.

Léo Gerville-Réache

University of Bordeaux - IMB UMR 52-51

**Abstract :** We analyze the main arguments that attempt to explain why there is no point in changing the envelope. Most people confuse estimation and calculation, conditional and unconditional probabilities, random and non-random variables, modelling and reality, chance and uncertainty (or credence), whatever *X* known and whatever *X* unknown, personal and interpersonal estimates. Such confusions suggest that there would be no point in swapping the envelope. Our analysis is that if the first envelope is opened, whatever the amount, it is rational and consistent to think that it is worth changing it. If the first envelope is not opened, it is rational and consistent to believe there is no point in swapping the envelope.

### Introduction

The two envelopes paradox is probably one of the most delicate in decision theory in recent decades. Many authors have addressed the issue and tried to explain the reasons for the problem in different ways. Currently, there is no consensus on a demonstration, since most people generally reject each other's demonstrations. Perhaps a consensus cannot be struck for the simple reason that one has rational interest in swapping the envelope. This article re-examines the original article in Wikipedia (2014-February), which is particularly exhaustive on the subject. Using footnotes, we try to demonstrate the various problems of interpretation and modelling.

************************************************************************

The two envelopes problem, also known as the exchange paradox, is a brain teaser, puzzle, or paradox in logic, philosophy, probability, and recreational mathematics. It is of special interest in decision theory, and for the Bayesian interpretation of probability theory. Historically, it arose as a variant of the necktie paradox.

### The problem:

You have two indistinguishable envelopes that each contain money. One contains twice as much as the other. You may pick one envelope and keep the money it contains. You pick at random, but before you open the envelope, you are offered the chance to take the other envelope instead.

It can be argued that it is to your advantage to swap envelopes by showing that your expected return on swapping exceeds the sum in your envelope. This leads to the paradoxical conclusion that it is beneficial to continue to swap envelopes indefinitely.

### Example

Assume the amount in my selected envelope is $20. If I happened to have selected the larger of the two envelopes ("larger" meaning the one with the larger amount of money), that would mean that the amount in my envelope is twice the amount in the other envelope. So in this case the amount in the other envelope would be $10 [1].

---

[1] Ok

However if I happened to have selected the smaller of the two envelopes, that would mean that the amount in the other envelope is twice the amount in my envelope. So in this second scenario the amount in the other envelope would be $40 [2].

The probability of either of these scenarios is one half, since there is a 50% chance that I initially happened to select the larger envelope and a 50% chance that I initially happened to select the smaller envelope. The expected value calculation for how much money is in the other envelope would be the amount in the first scenario times the probability of the first scenario plus the amount in the second scenario times the probability of the second scenario, which is $10*1/2 + $40*1/2.[3] The result of this calculation is that the expected value [4] of money in the other envelope is $25. Since this is greater than my selected envelope, it would appear to my advantage to always switch envelopes [5].

A large number of solutions have been proposed. The usual scenario is that one writer proposes a solution that solves the problem as stated, but then another writer discovers that altering the problem slightly revives the paradox. In this way, a family of closely related formulations of the problem have been created, which are discussed in the literature.

No proposed solution is widely accepted as correct.[1] Despite this it is common for authors to claim that the solution to the problem is easy, even elementary.[2] However, when investigating these elementary solutions they often differ from one author to the next. Since 1987 new papers have been published every year.[3] [6]

**Problem**

Basic setup: You are given two indistinguishable envelopes, each of which contains a positive sum of money. One envelope contains twice as much as the other. You may pick one envelope and keep whatever amount it contains. You pick one envelope at random but before you open it you are given the chance to take the other envelope instead.[4]

The switching argument: Now suppose you reason as follows:

1) I denote by *A* the amount in my selected envelope. [7]

2) The probability that *A* is the smaller amount is 1/2, and that it is the larger amount is also 1/2. [8]

3) The other envelope may contain either 2*A* or *A*/2. [9]

4) If *A* is the smaller amount, then the other envelope contains 2*A*. [10]

---

[2] Ok

[3] This is an estimate of the expected gain, provided the opened envelope contains $20.

[4] Not "the expected value" but an estimation of the expected value.

[5] If my utility is linear, Ok.

[6] It is pleasing to see no consensus to explain that the reasoning based on the estimated expectation, which is conditional on the value of the envelope selected, is wrong. Indeed, there is nothing wrong or inconsistent in this conditional estimate. The error in the above reasoning is to think that one may calculate an expectation, whereas it is an estimate of the conditional expectation in swapping the envelope.

[7] Here *A* is unknown! One should not confuse an analysis that would be consistent whatever *A* known with one whatever *A* unknown.

[8] Ok.

[9] Ok but these values are unknown.

[10] Ok.

5) If *A* is the larger amount, then the other envelope contains *A*/2. [11]

6) Thus the other envelope contains 2*A* with probability 1/2 and *A*/2 with probability 1/2. [12]

7) So the expected value of the money in the other envelope is

$$1\backslash 2*2A + 1\backslash 2*A\backslash 2 = 5\backslash 4*A \quad [13]$$

1) This is greater than *A*, [14] so I gain on average by swapping. [15]

2) After the switch, I can denote that content by *B* and reason in exactly the same manner as above. [16]

3) I will conclude that the most rational thing to do is to swap back again. [17]

4) To be rational, I will thus end up swapping envelopes indefinitely. [18]

5) As it seems more rational to open just any envelope than to swap indefinitely, we have a contradiction. [19]

The puzzle: The puzzle is to find the flaw in the very compelling line of reasoning above.

**Common resolution**

A common way to resolve the paradox, both in popular literature and in the academic literature in philosophy, is to observe that A stands for different things at different places in the expected value calculation, step 7 above.[5] In the first term A is the smaller amount while in the second term A is the larger amount. To mix different instances of a variable in the same formula like this is said to be illegitimate [20], so step 7 is incorrect, and this is the cause of the paradox.

According to this analysis, a correct argument runs on the following lines. By definition the amount in one envelope is twice as much as in the other. Denoting the lower of the two amounts by *X*, we write the expected value calculation as

$$1\backslash 2*2X + 1\backslash 2*X = 3\backslash 2*X \quad [21]$$

---

[11] Ok.

[12] Ok, but it is a modelling (or an estimate) of probabilities of scenarios.

[13] This is an estimate of expectation whose value is unknown!

[14] This is not true. As A is unknown, one cannot say that 5/4*A is greater than A. One can only say that 5/4*A>A, whatever A known, A being greater than 0. This distinction is fundamental.

[15] This is not true. If *A* is unknown, it is by this argument that some authors speak about infinite expectation for *A*. One cannot compare *A* and 5/4*A. It is also necessary to define precisely the random experiment, which will be repeated, to be able to talk about results "on average".

[16] I can do it but it would be equally wrong.

[17] Inexact.

[18] Inexact.

[19] Reasoning without knowledge of *A*, we also obtain infinite expectation…

[20] Unless we know *A*.

[21] Ok, this is the expectation of the game before randomly selecting an envelope. However, this expectation is unknown because *X* is unknown and the whole issue of the paradox is to understand how one should estimate the expectation, conditionally (or not) on the fact that the amount of the first envelope is known.

Here *X* stands for the same thing in every term of the equation. We learn that 1.5*X*, is the average expected value in either of the envelopes. Being less than 2*X*, the value in the greater envelope, there is no reason to swap the envelopes. [22]

**Mathematical details**

Let us rewrite the preceding calculations in a more detailed notation that explicitly distinguishes random from not-random quantities (a different distinction from the usual in ordinary, deterministic mathematics between variables and constants [23]. This is a useful way to compare with the next, alternative, resolution. So far, we were thinking of the two amounts of money in the two envelopes as fixed [24]. The only randomness lies in which one goes into which envelope. [25] We called the smaller amount *X*, let us denote the larger amount by *Y*. Given the values x and y of *X* and *Y*, where *y* = 2*x* and *x* > 0, the problem description tells us (whether or not *x* and *y* are known)

$$P(A=x|X=x) = P(A=y|X=x) = 1\backslash 2 \quad [26]$$

for all possible values *x* of the smaller amount *X*; there is a corresponding definition of the probability distribution of *B* given *X* and *Y*. In our resolution of the paradox, we guessed that in Step 7 the writer was trying to compute the expected value of *B* given *X*=*x*. Splitting the calculation over the two possibilities for which envelope contains the smaller amount, it is certainly correct to write

$$E(B|X=x) = E(B|B>A, X=x) * P(B>A|X=x) + E(B|B<A, X=x) * P(B<A|X=x). \quad [27]$$

At this point, the writer correctly substitutes the value 1/2 for both of the conditional probabilities on the right hand side of this equation (Step 2). At the same time he correctly substitutes the random variable *B* inside the first conditional expectation for 2*A*, when taking its expectation value given *B*>*A* and *X* = *x*, and he similarly correctly substitutes the random variable *B* for *A*/2 when taking its expectation value given *B*<*A* and *X* = *x* (Steps 4 and 5). He would then arrive at the completely correct equation

$$E(B|X=x) = 1/2*E(2A|B>A, X=x) + 1/2*E(A/2|B<A, X=x)$$

$$= E(A|B>A, X=x) + 1/4* E(A|B<A, X=x)$$

However he now proceeds, in the first of the two terms on the right hand side, to replace the expectation value of *A* given that Envelope A contains the smaller amount and given that the amounts are *x* and 2*x*, by the random quantity *A* itself. Similarly, in the second term on the right hand side he replaces the expectation value of *A* given now that Envelope A contains the larger amount and given that the amounts are *x* and 2*x*, also by the random quantity *A* itself. The correct substitutions would have been, of course, *x* and 2*x* respectively, leading to a correct conclusion

---

[22] This expectation tells us that as long as we cannot estimate this expectation, there is no reason to prefer one or the other envelopes. However, when we know the amount in the first envelope, we would prefer the other and that would be perfectly rational…

[23] This is an interpretation of the difference between "chance" and "uncertainty" which is not scientifically or philosophically resolved.

[24] It is clear that when an envelope is selected, the amounts exist and are fixed but are unknown.

[25] What do we do about the uncertainty of the amounts in the envelopes?

[26] This is not clear; x is known or not? A conditional probability is a number and therefore x must be known. Otherwise, the probability P(*A*=x|*X*=x) is the following random variable : $\mathbf{1}_{\{X=A\}}$.

[27] This is not clear; we have random variables or numbers? What can be written is that if you opened the first envelope, either we chose *A* and then E(*B*|*A*=*a*) = 3/2*\**a*, or we chose *B* and then E(*A*|*B*=*b*) = 3/2*\**b*. These are the two possible estimates. As we know that max(*a*, *b*) = 2*\*min(*a*, *b*), we necessarily have E(*B*|*A*=*a*) ≠ E(*A*|*B*=*b*). There is no problem here.

$$E(B|X=x) = 3/2*x.\ ^{28}$$

Naturally this coincides with the expectation value of *A* given *X=x*. [29]

Indeed, in the two contexts in which the random variable A appears on the right hand side, it is standing for two different things, since its distribution has been conditioned on different events. Obviously, A tends to be larger, when we know that it is greater than B and when the two amounts are fixed, and it tends to be smaller, when we know that it is smaller than B and the two amounts are fixed, cf. Schwitzgebel and Dever (2007, 2008). In fact, it is exactly twice as large in the first situation as in the second situation.

The preceding resolution was first noted by Bruss in 1996.[6] Falk gave a concise exposition in 2009.[7] [30]

**Alternative interpretation**

The first solution above does not explain what is wrong if the player is allowed to open the first envelope before offered the option to switch [31]. In this case, *A* stands for the value that then seen throughout all subsequent calculations. The mathematical variable A stands for any particular amount he might see there (it is a mathematical variable, a generic possible value of a random variable). The reasoning appears to show that whatever amount he would see there, he would decide to switch. [32]

Hence, he does not need to look in the envelope at all: he knows that if he looks, and goes through the calculations, they will tell him to switch, whatever he saw in the envelope. [33]

In this case, at Steps 6, 7 and 8 of the reasoning, A is any fixed possible value of the amount of money in the first envelope.

Thus, the proposed "common resolution" above breaks down, and another explanation is needed.

This interpretation of the two envelopes problem appears in the first publications in which the paradox was introduced, Gardner (1989) and Nalebuff (1989). It is common in the more mathematical literature on the problem.

The "common resolution" above depends on a particular interpretation of what the writer of the argument is trying to calculate: namely, it assumes he is after the (unconditional) expectation value of what's in Envelope B. [34] In the mathematical literature on Two Envelopes Problem (and in particular, in the literature where it was first introduced to the world), another interpretation is more common, involving the conditional expectation value (conditional on what might be in Envelope A). To solve this and related interpretations or versions of the problem, most authors use the Bayesian interpretation of probability. [35]

**Introduction to resolutions based on Bayesian probability theory**

Here the ways the paradox can be resolved depend to a large degree on the assumptions that are made about the things that are not made clear in the setup and the proposed argument for switching.[8][9] The most usual assumption about the way the envelopes are set up is that a sum of money is in one envelope, and twice that sum is in another envelope. One of the two envelopes is randomly given to the player (envelope *A*). It is not made clear exactly how the smaller of the two sums is determined, what values it could possibly take and, in particular, whether there is a maximum sum it might contain. [36] [10][11]. It is also not specified whether the player can look in Envelope A before deciding whether or not to switch. A further ambiguity in the paradox is that it is not made clear in the proposed argument whether the amount A in Envelope A is intended to be a constant, a random variable, or some other quantity. [37]

---

[28] Here, *x* is unknown and is therefore modelled as a random variable but with unknown distribution ... If *x* is known here, the game has no interest because we know the smaller amount!

[29] (Idem...)

[30] The attempted demonstration is not admissible.

[31] The first "solution" explains nothing.

[32] This is the fundamental difference between "whatever x known" and "whatever x unknown".

[33] By seeing the amount of the first envelope, we are able to calculate estimations and that is what is essential, even if the decision does not depend on the result of this estimate.

[34] Unconditional expectation (in the sense that the amount of the first envelope is unknown) is a random variable with unknown distribution!

[35] Bayesian?

If it assumed that there is a maximum sum that can be put in the first envelope, then a simple and mathematically sound resolution is possible within the second interpretation [38]. Step 6 in the proposed line of reasoning is not always true, since if the player holds more than the maximum sum that can be put into the first envelope they must hold the envelope containing the larger sum, and are thus certain to lose by switching. This may not occur often, but when it does, the heavy loss the player incurs means that, on average, there is no advantage in switching. This resolves all practical cases of the problem, whether or not the player looks in their envelope [12] [39].

It can be envisaged, however, that the sums in the two envelopes are not limited. [40] This requires a more careful mathematical analysis, and also uncovers other possible interpretations of the problem. If, for example, the smaller of the two sums of money is considered equally likely to be one of infinitely many positive integers, without upper limit—the probability that it is any given number is always zero. This absurd situation exemplifies an improper prior, and this is generally considered to resolve the paradox in this case. [41]

It is possible to devise a distribution for the sums possible in the first envelope such that the maximum value is unlimited, computation of the expectation of what is in B given what is in A seems to dictate you should switch, and the distribution constitutes a proper prior.[13] In these cases it can be shown that the expected sum in both envelopes is infinite. There is no gain, on average, in swapping. [42]

The first two resolutions we present correspond, technically speaking, first to A being a random variable, and secondly to it being a possible value of a random variable (and the expectation being computed is a conditional expectation). At the same time, in the first resolution the two original amounts of money seem to be thought of as being fixed, while in the second they are also thought of as varying. Thus there are two main interpretations of the problem, and two main resolutions. [43]

**Proposed resolutions to the alternative interpretation**

Nalebuff (1989), Christensen and Utts (1992), Falk and Konold (1992), Blachman, Christensen and Utts (1996),[14] Nickerson and Falk (2006), pointed out that if the amounts of money in the two envelopes have any proper probability distribution representing the player's prior beliefs about the amounts of money in the two envelopes, then it is impossible that whatever the amount A=a in the first envelope might be, it would be equally likely, according to these prior beliefs, that the second contains a/2 or 2a. Thus step 6 of the argument, which leads to always switching, is a non-sequitur. [44]

---

[36] It is a theory about the probability laws of sums of money at stake in the problem of two envelopes. Where is the a-priori?

[37] The description is very partial against this approach, even if I agree with it.

[38] Ok

[39] Here, the decision to swap depends on the a-prior.

[40] Infinity is poking its nose out with all the paradoxes it may generate. More is needed on this question…

[41] Let's be reasonable, there is a finite amount in each envelope. It is possible!

[42] Fortunately this is not the solution of the two envelopes paradox.

[43] The attempted demonstration is not admissible.

*Mathematical details*

According to this interpretation, the writer is carrying out the following computation, where he is conditioning now on the value of A, the amount in Envelope A, not on the pair amounts in the two envelopes X and Y:

$$E(B|A=a) = E(B|A=a, A<B) * P(A<B|A=a) + E(B|A=a, A>B) * P(A>B|A=a).$$ [45]

Completely correctly, and according to Step 5, the two conditional expectation values are evaluated [46] as

$$E(B|A=a, A<B) = E(2A|A=a, A<B) = 2a,$$

$$E(B|A=a, A>B) = E(A/2|A=a, A>B) = a/2.$$

However in Step 6 the writer is invoking Steps 2 and 3 to get the two conditional probabilities, and effectively replacing the two conditional probabilities of Envelope *A* containing the smaller and larger amount, respectively, given the amount actually in that envelope, both by the unconditional probability 1/2: he makes the substitutions

$$P(A<B|A=a) = 1/2,$$ [47]

$$P(A>B|A=a) = 1/2.$$ [48]

But intuitively, we would expect that the larger the amount in A, the more likely it is the larger of the two, and vice-versa. And it is a mathematical fact, as we will see in a moment, that it is impossible that both of these conditional probabilities are equal to 1/2 for all possible values of a [49]. In fact, for step 6 to be true, whatever a might be, the smaller amount of money in the two envelopes must be equally likely to be between 1 and 2, as between 2 and 4, as between 4 and 8, ... ad infinitum [50]. But there is no way to divide total probability 1 into an infinite number of pieces that are not only all equal to one another, but also all larger than zero [51]. Yet the smaller amount of money in the two envelopes must have probability larger than zero to be in at least one of the just mentioned ranges.

To see this, suppose that the chance that the smaller of the two envelopes contains an amount between $2^n$ and $2^{n+1}$ is $p(n)$, where n is any whole number, positive or negative, and for definiteness we include the lower limit but exclude the upper in each interval. It follows that the conditional probability that the envelope in our hands contains the smaller amount of money of the two, given that its contents are between $2^n$ and $2^{n+1}$, is

$$(p(n)/2) / (p(n)/2 + p(n-1)/2).$$

If this is equal to 1/2, it follows by simple algebra that

$$2p(n)=p(n)+p(n-1),$$

---

[44] "Whatever the amount *A=a*" does not mean that all amounts are possible but, whatever the amount discovered, I do not have to probabilize the possible amounts of the first envelope but have to work conditionally on this known amount...

[45] Ok, this is the estimate of the expectation conditionally on the knowledge of amount *a* in envelope *A*.

[46] More precisely: estimated.

[47] Ok, this is an estimate.

[48] Ok, this is an estimate.

[49] How do we show it? To be continued...

[50] Infinity again!

[51] Ok, but that is not the problem.

or $p(n)=p(n-1)$. This must be true for all *n*, an impossibility. [52]

**New variant**

Though Bayesian probability theory can resolve the alternative interpretation of the paradox above, it turns out that examples can be found of proper probability distributions, such that the expected value of the amount in the second envelope given that in the first does exceed the amount in the first, whatever it might be. The first such example was already given by Nalebuff (1989). See also Christensen and Utts (1992)[15][16][17][18]

Denote again the amount of money in the first envelope by *A* and that in the second by *B*. We think of these as random. Let *X* be the smaller of the two amounts and *Y*=2*X* be the larger. Notice that once we have fixed a probability distribution for *X* then the joint probability distribution of *A*,*B* is fixed, since *A*,*B* = *X*,*Y* or *Y*,*X* each with probability 1/2, independently of *X*,*Y*.

The bad step 6 in the "always switching" argument led us to the finding $E(B|A=a)>a$ for all *a*, and hence to the recommendation to switch, whether or not we know *a* [53]. Now, it turns out that one can quite easily invent proper probability distributions for *X*, the smaller of the two amounts of money, such that this bad conclusion is still true. One example is analysed in more detail, in a moment.

It cannot be true that whatever *a*, given *A*=*a*, *B* is equally likely to be *a*/2 or 2*a*, but it can be true that whatever *a*, given *A*=*a*, *B* is larger in expected value than *a*. [54]

Suppose for example (Broome, 1995)[19] that the envelope with the smaller amount actually contains 2n dollars with probability $2n/3n+1$ where n = 0, 1, 2,... [55] These probabilities sum to 1, hence the distribution is a proper prior (for subjectivists) and a completely decent probability law also for frequentists.

Imagine what might be in the first envelope. A sensible strategy would certainly be to swap when the first envelope contains 1, as the other must then contain 2. Suppose on the other hand the first envelope contains 2. In that case there are two possibilities: the envelope pair in front of us is either {1, 2} or {2, 4}. All other pairs are impossible. The conditional probability that we are dealing with the {1, 2} pair, given that the first envelope contains 2, is

$$P(\{1,2\})\backslash 2 = (P(\{1,2\})/2) / (\{P(\{1,2\})/2 + P(\{2,4\})/2\}$$
$$= P(\{1,2\}) / (\{P(\{1,2\})+P(\{2,4\})\})$$
$$= (1/3) / (1/3 + 2/9) = 3/5,$$

and consequently the probability it's the {2, 4} pair is 2/5, since these are the only two possibilities. In this derivation, P({1,2})/2 is the probability that the envelope pair is the pair 1 and 2, and Envelope *A* happens to contain 2; P({2,4})/2 is the probability that the envelope pair is the pair 2 and 4, and (again) Envelope *A* happens to contain 2. Those are the only two ways that Envelope *A* can end up containing the amount 2.

It turns out that these proportions hold in general unless the first envelope contains 1. Denote by a the amount we imagine finding in Envelope *A*, if we were to open that envelope, and suppose that $a=2n$ for some n ≥ 1. In that case the other envelope contains *a*/2 with probability 3/5 and 2*a* with probability 2/5.

So either the first envelope contains 1, in which case the conditional expected amount in the other envelope is 2, or the first envelope contains *a* > 1, and though the second envelope is more likely to be smaller than larger, its conditionally expected amount is larger: the conditionally expected amount in Envelope *B* is

$$3/5*a/2 + 2/5*2a = 11/10*a$$

which is more than a. This means that the player who looks in Envelope *A* would decide to switch whatever he saw there. Hence there is no need to look in Envelope *A* to make that decision.

---

[52] The attempted demonstration is not admissible.

[53] If you don't know *a*, E(*B*|*A* = *a*) is a random variable and the mathematical expression is not clear; if we know *a*, E(*B*|*A* = *a*) is an estimate, its value is indeed 5/4*\**a*.

[54] No, it is true because it's an estimate.

[55] At this rate, I can also assume that the smallest envelope contains 100$. Then, if I first select 100$, I swap and if I select first 200$, I keep it.

This conclusion is just as clearly wrong as it was in the preceding interpretations of the Two Envelopes Problem. But now the flaws noted above do not apply; the *a* in the expected value calculation is a constant and the conditional probabilities in the formula are obtained from a specified and proper prior distribution. [56]

**Proposed resolutions**

Some writers think that the new paradox can be defused.[20] Suppose E(*B*|*A*=*a*)>*a* for all *a*. As remarked before, this is possible for some probability distributions of *X* (the smaller amount of money in the two envelopes). Averaging over a, it follows either that E(*B*)>E(*A*), or alternatively that E(*B*)=E(*A*)=∞. But *A* and *B* have the same probability distribution, and hence the same expectation value, by symmetry (each envelope is equally likely to be the smaller of the two). Thus both have infinite expectation values, and hence so must *X* too.

Thus if we switch for the second envelope because its conditional expected value is larger than what actually is in the first, whatever that might be, we are exchanging an unknown amount of money whose expectation value is infinite for another unknown amount of money with the same distribution and the same infinite expected value [57]. The average amount of money in both envelopes is infinite. Exchanging one for the other simply exchanges an average of infinity with an average of infinity [58].

Probability theory therefore tells us why and when the paradox can occur and explains to us where the sequence of apparently logical steps breaks down. In this situation, Steps 6 and Steps 7 of the standard Two Envelopes argument can be replaced by correct calculations of the conditional probabilities that the other envelope contains half or twice what's in A, and a correct calculation of the conditional expectation of what's in B given what's in A. Indeed, that conditional expected value is larger than what's in A. But because the unconditional expected amount in A is infinite, this does not provide a reason to switch, because it does not guarantee that on average you'll be better off after switching. One only has this mathematical guarantee in the situation that the unconditional expectation value of what's in A is finite. But then the reason for switching without looking in the envelope, E(*B*|*A*=*a*)>*a* for all *a*, simply cannot arise [59].

Many economists prefer to argue that in a real-life situation, the expectation of the amount of money in an envelope cannot be infinity, for instance, because the total amount of money in the world is bounded; therefore any probability distribution describing the real world would have to assign probability 0 to the amount being larger than the total amount of money on the world [60]. Therefore the expectation of the amount of money under this distribution cannot be infinity. The resolution of the second paradox, for such writers, is that the postulated probability distributions cannot arise in a real-life situation. These are similar arguments as used to explain the St. Petersburg Paradox [61] [62].

---

[56] The attempted demonstration is not admissible.

[57] Yes, this means that if you do not know the amount of the first envelope, reasoning with random variables means choosing between two infinite expectations.

[58] Ok.

[59] Be careful, if *a* is known, then, "whatever *a* is", the expectation is finite! That is why there is a fundamental difference between "whatever *a* known" and "whatever *a* unknown".

[60] There is no economic question in the two envelopes paradox.

[61] The paradox of St Petersburg discusses the relevance of the expected gain as a criterion for "rational" decision. It is very interesting but it is not the main problem here.

[62] The attempted demonstration is admissible, when you don't know the amount of the first envelope.

**Foundations of mathematical economics**

In mathematical economics and the theory of utility [63], which explains economic behaviour in terms of expected utility, there remains a problem to be resolved.[21] In the real world we presumably would not indefinitely exchange one envelope for the other (and probability theory, as just discussed, explains quite well why calculations of conditional expectations might mislead us). Yet the expected utility based theory of economic behaviour assumes that people do (or should) make economic decisions by maximizing expected utility, conditional on present knowledge. If the utility function is unbounded above [64], then the theory can still predict infinite switching.

Fortunately for mathematical economics and the theory of utility, it is generally agreed that as an amount of money increases, its utility to the owner increases less and less [65], and ultimately there is a finite upper bound to the utility of all possible amounts of money. We can pretend that the amount of money in the whole world is as large as we like, yet the utility that the owner of all that money experiences, while rising further and further, will never rise beyond a certain point no matter how much is in his possession. For decision theory and utility theory, the two envelope paradox illustrates that unbounded utility does not exist in the real world, so fortunately there is no need to build a decision theory that allows unbounded utility, let alone utility of infinite expectation.[citation needed] [66].

**Controversy among philosophers**

As mentioned above, any distribution producing this variant of the paradox must have an infinite mean. So before the player opens an envelope the expected gain from switching is "$\infty - \infty$", [67] which is not defined. In the words of Chalmers this is "just another example of a familiar phenomenon, the strange behaviour of infinity" [68].[22] Chalmers suggests that decision theory generally breaks down when confronted with games having a diverging expectation, and compares it with the situation generated by the classical St. Petersburg paradox [69].

However, Clark and Shackel argue that this blaming it all on "the strange behaviour of infinity" does not resolve the paradox at all; neither in the single case nor the averaged case. They provide a simple example of a pair of random variables both having infinite mean but where it is clearly sensible to prefer one to the other, both conditionally and on average.[23] They argue that decision theory should be extended so as to allow infinite expectation values in some situations. [70] [71]

---

[63] Forget the utility or assume that it is linear in the two envelopes problem. The question (very interesting in fact) is not to explain what people would do but how they have to reason. Everyone is free to decide later depending on his personal utility.

[64] Still it's curious to cross the infinite in a science that must stick to modeling reality.

[65] The logarithm still tends to infinity …

[66] The attempted demonstration is not admissible.

[67] This is actually a possible modelling of ignorance at this step of the game.

[68] Rather it is the best mathematical modelling of our ignorance when the first envelope is not open. This translation simply tells us that swapping before knowing the amount of the first envelope is unfounded.

[69] It's another problem.

[70] I agree, it's interesting to imagine that although two infinite sums are mathematically equal, we prefer, as humans, the amount multiplied by 5/4.

[71] The attempted demonstration is not admissible but this view of the ambiguous relationship between human beings and infinity is interesting.

**Non-probabilistic variant**

The logician Raymond Smullyan questioned if the paradox has anything to do with probabilities at all. [72] [24] He did this by expressing the problem in a way that does not involve probabilities. [73] The following plainly logical arguments lead to conflicting conclusions:

Let the amount in the envelope chosen by the player be A. By swapping, the player may gain *A* or lose *A*/2. So the potential gain is strictly greater than the potential loss.

Let the amounts in the envelopes be *X* and 2*X*. Now by swapping, the player may gain *X* or lose *X*. So the potential gain is equal to the potential loss. [74]

**Proposed resolutions**

A number of solutions have been put forward. Careful analyses have been made by some logicians. Though solutions differ, they all pinpoint semantic issues concerned with counterfactual reasoning. We want to compare the amount that we would gain by switching if we would gain by switching, with the amount we would lose by switching if we would indeed lose by switching. However, we cannot both gain and lose by switching at the same time [75]. It should not be confused chance and uncertainty... We are asked to compare two incompatible situations. Only one of them can factually occur, the other is a counterfactual situation- somehow imaginary. To compare them at all, we must somehow "align" the two situations, providing some definite points in common.

James Chase (2002) argues that the second argument is correct because it does correspond to the way to align two situations (one in which we gain, the other in which we lose), which is preferably indicated by the problem description.[25] Also Bernard Katz and Doris Olin (2007) argue this point of view.[26] In the second argument, we consider the amounts of money in the two envelopes as being fixed; what varies is which one is first given to the player. Because that was an arbitrary and physical choice, the counterfactual world in which the player, counterfactually, got the other envelope to the one he was actually (factually) given is a highly meaningful counterfactual world and hence the comparison between gains and losses in the two worlds is meaningful. This comparison is uniquely indicated by the problem description, in which two amounts of money are put in the two envelopes first, and only after that is one chosen arbitrarily and given to the player [76]. In the first argument, however, we consider the amount of money in the envelope first given to the player as fixed and consider the situations where the second envelope contains either half or twice that amount [77]. This would only be a reasonable counterfactual world if in reality the envelopes had been filled as follows: first, some amount of money is placed in the specific envelope that will be given to the player; and secondly, by some arbitrary process, the other envelope is filled (arbitrarily or randomly) either with double or with half of that amount of money [78]. Byeong-Uk Yi (2009), on the other hand, argues that comparing the amount you would gain if you would gain by switching with the amount you would lose if you would lose by switching is a meaningless exercise from the outset.[27] According to his analysis, all three implications (switch, indifferent, do not switch) are incorrect [79]. He analyses Smullyan's arguments in detail, showing that intermediate steps are being taken, and pinpointing exactly where an incorrect inference is made according to his formalization of counterfactual inference. An important difference with Chase's analysis is that he does not take account of the part of the story where we are told that the envelope called Envelope *A* is decided completely at random. Thus, Chase puts probability back into the problem description in order to conclude that arguments 1 and 3 are incorrect, argument 2 is correct, while Yi keeps "two envelope problem without probability" completely free of probability, and comes to the conclusion that there are no reasons to prefer any action. This corresponds to the view of Albers et al., that without probability ingredient, there is no way to argue that one action is better than another, anyway.

---

[72] Be careful, what does "probabilities" mean? If they only measure "chance", the question arises. But here we are dealing with uncertainty and decision.

[73] Probability measures not only chance, and so much the better.

[74] The attempted demonstration is not admissible.

[75] Ok, and a coin cannot fall on both Heads and Tails at the same time but that does not stop us saying P(Heads) = P(Tails) = 1/2. We mustn't confuse chance and uncertainty.

[76] It is indeed unsettling to note that the value of the estimate of the conditional expectation of the two envelopes problem is equal to the value of the conditional expectation of a game where the values of the envelope not selected would be randomly assigned after selecting and opening the first envelope. However, in the first case, it is an estimate and in the second, it is a probabilistic calculation.

[77] It is a rational and consistent modelling of latent uncertainty.

In perhaps the most recent paper on the subject, Bliss argues that the source of the paradox is that when one mistakenly believes in the possibility of a larger payoff that does not, in actuality, exist, one is mistaken by a larger margin than when one believes in the possibility of a smaller payoff that does not actually exist.[28] If, for example, the envelopes contained $5.00 and $10.00 respectively, a player who opened the $10.00 envelope would expect the possibility of a $20.00 payout that simply does not exist [80]. Were that player to open the $5.00 envelope instead, he would believe in the possibility of a $2.50 payout, which constitutes a smaller deviation from the true value.

Albers, Kooi, and Schaafsma (2005) consider that without adding probability (or other) ingredients to the problem, Smullyan's arguments do not give any reason to swap or not to swap, in any case. Thus there is no paradox. This dismissive attitude is common among writers from probability and economics: Smullyan's paradox arises precisely because he takes no account whatever of probability or utility [81].

**Extensions to the problem**

Since the two envelopes problem became popular, many authors have studied the problem in depth in the situation in which the player has a prior probability distribution of the values in the two envelopes, and does look in Envelope A. One of the most recent such publications is by McDonnell and Douglas (2009), who also consider some further generalizations.[29]

If a priori we know that the amount in the smaller envelope is a whole number of some currency units, then the problem is determined, as far as probability theory is concerned, by the probability mass function $p(x)$ describing our prior beliefs that the smaller amount is any number x = 1,2, ... ; the summation over all values of *x* being equal to 1. It follows that given the amount *a* in Envelope *A*, the amount in Envelope *B* is certainly 2*a* if a is an odd number. However, if *a* is even, then the amount in Envelope *B* is 2*a* with probability $p(a)/(p(a/2)+p(a))$, and *a*/2 with probability $p(a/2)/(p(a/2)+p(a))$. If one would like to switch envelopes if the expectation value of what is in the other is larger than what we have in ours, then a simple calculation shows that one should switch if $p(a/2) < 2p(a)$, keep to Envelope *A* if $p(a/2) > 2p(a)$.

If on the other hand the smaller amount of money can vary continuously, and we represent our prior beliefs about it with a probability density *f(x)*, thus a function that integrates to one when we integrate over *x* running from zero to infinity, then given the amount a in Envelope *A*, the other envelope contains 2*a* with probability $2f(a)/(f(a/2)+2f(a))$, and *a*/2 with probability $f(a/2)/(f(a/2)+2f(a))$. If again we decide to switch or not according to the expectation value of what's in the other envelope, the criterion for switching now becomes $f(a/2) < 4f(a)$.

---

[78] We still confuse "chance" and "uncertainty" and it is said that one should not probabilize uncertainty.

[79] So, it would be a paradox because adding the fact that one is twice the other, we can do nothing more and if we do not know the rate between the two envelopes, there is no problem.

[80] This argument is a straight concept of uncertainty; it is a chronologist concept of chance.

[81] The attempted demonstration is not admissible.

The difference between the results for discrete and continuous variables may surprise many readers. Speaking intuitively, this is explained as follows. Let $h$ be a small quantity and imagine that the amount of money we see when we look in Envelope *A* is rounded off in such a way that differences smaller than $h$ are not noticeable, even though actually it varies continuously. The probability that the smaller amount of money is in an interval around *a* of length *h*, and Envelope *A* contains the smaller amount is approximately $f(a)*h*(1/2)$. The probability that the larger amount of money is in an interval around *a* of length h corresponds to the smaller amount being in an interval of length $h/2$ around $a/2$. Hence the probability that the larger amount of money is in a small interval around *a* of length h and Envelope *A* contains the larger amount is approximately $f(a/2)*(h/2)*(1/2)$. Thus, given Envelope *A* contains an amount about equal to *a*, the probability it is the smaller of the two is roughly $f(a)*h*(1/2) / (f(a)*h*(1/2) + f(a/2)*(h/2)*(1/2)) = 2f(a) / (2f(a) + f(a/2))$.

If the player only wants to end up with the larger amount of money, and does not care about expected amounts, then in the discrete case he should switch if *a* is an odd number, or if a is even and $p(a/2) < p(a)$. In the continuous case he should switch if $f(a/2) < 2f(a)$. [82]

Some authors prefer to think of probability in a frequentist sense. If the player knows the probability distribution used by the organizer to determine the smaller of the two values, then the analysis would proceed just as in the case when p or f represents subjective prior beliefs (Let). However, what if we take a frequentist point of view, but the player does not know what probability distribution is used by the organiser to fix the amounts of money in any one instance? Thinking of the arranger of the game and the player as two parties in a two person game, puts the problem into the range of game theory. The arranger's strategy consists of a choice of a probability distribution of *x*, the smaller of the two amounts. Allowing the player also to use randomness in making his decision, his strategy is determined by his choosing a probability of switching $q(a)$ for each possible amount of money a he might see in Envelope *A*.[83] In this section we so far only discussed fixed strategies, that is strategies for which q only takes the values 0 and 1, and we saw that the player is fine with a fixed strategy, if he knows the strategy of the organizer. In the next section we will see that randomized strategies can be useful when the organizer's strategy is not known. [84]

**Randomized solutions**

Suppose as in the previous section that the player is allowed to look in the first envelope before deciding whether to switch or to stay. [85] We'll think of the contents of the two envelopes as being two positive numbers, not necessarily two amounts of money. The player is allowed either to keep the number in Envelope A, or to switch and take the number in Envelope B. We'll drop the assumption that one number is exactly twice the other, we'll just suppose that they are different and positive. On the other hand, instead of trying to maximize expectation values, we'll just try to maximize the chance that we end up with the larger number. [86]

In this section we ask the question, is it possible for the player to make his choice in such a way that he goes home with the larger number with probability strictly greater than half, however the organizer has filled the two envelopes?

---

[82] It's another problem...

[83] The attempted demonstration is not admissible.

[84] We shall see...

[85] Ok.

[86] Why? In this case E(B|A=x)=x, and therefore changing doesn't increase the conditional expectation of gain.

We are given no information at all about the two numbers in the two envelopes, except that they are different, and strictly greater than zero. The numbers were written down on slips of paper by the organiser, put into the two envelopes. The envelopes were then shuffled, the player picks one, calls it Envelope *A*, and opens it.[87]

We are not told any joint probability distribution of the two numbers. We are not asking for a subjectivist solution. We must think of the two numbers in the envelopes as chosen by the arranger of the game according to some possibly random procedure, completely unknown to us, and fixed. Think of each envelope as simply containing a positive number and such that the two numbers are not the same. The job of the player is to end up with the envelope with the larger number. This variant of the problem, as well as its solution, is attributed by McDonnell and Abbott, and by earlier authors, to information theorist Thomas M. Cover.[30]

Counter-intuitive though it might seem, there is a way that the player can decide whether to switch or to stay so that he has a larger chance than 1/2 of finishing with the bigger number [88], however the two numbers are chosen by the arranger of the game. However, it is only possible with a so-called randomized algorithm: the player must be able to generate his own random numbers. Suppose he is able to produce a random number, let's call it *Z*, such that the probability that *Z* is larger than any particular quantity *z* is exp(-*z*). Note that exp(-*z*) starts off equal to 1 at *z*=0 and decreases strictly and continuously as z increases, tending to zero as *z* tends to infinity. So the chance is 0 that *Z* is exactly equal to any particular number, and there is a positive probability that *Z* lies between any two particular different numbers [89]. The player compares his *Z* with the number in Envelope A. If *Z* is smaller he keeps the envelope. If *Z* is larger he switches to the other envelope.

Think of the two numbers in the envelopes as fixed (though of course unknown to the player). Think of the player's random *Z* as a probe with which he decides whether the number in Envelope A is small or large. If it is small compared to *Z* he switches, if it is large compared to *Z* he stays.

If both numbers are smaller than the player's *Z*, his strategy does not help him. He ends up with the Envelope B, which is equally likely to be the larger or the smaller of the two. If both numbers are larger than *Z* his strategy does not help him either, he ends up with the first Envelope *A*, which again is equally likely to be the larger or the smaller of the two. However if *Z* happens to be in between the two numbers, then his strategy leads him correctly to keep Envelope A if its contents are larger than those of *B*, but to switch to Envelope *B* if A has smaller contents than *B*. Altogether, this means that he ends up with the envelope with the larger number with probability strictly larger than 1/2. To be precise, the probability that he ends with the "winning envelope" is 1/2 + P(*Z* falls between the two numbers)/2. [90]

In practice, the number *Z* we have described could be determined to the necessary degree of accuracy as follows. Toss a fair coin many times, and convert the sequence of heads and tails into the binary representation of a number *U* between 0 and 1: for instance, *HTHHTH*... becomes the binary representation of *u*=0.101101.. . In this way, we generate a random number U, uniformly distributed between 0 and 1. Then define $Z = -\ln(U)$ where "ln" stands for natural logarithm, i.e., logarithm to base e. Note that we just need to toss the coin long enough to verify whether *Z* is smaller or larger than the number a in the first envelope—we do not need to go on for ever. We only need to toss the coin a finite (though random) number of times: at some point we can be sure that the outcomes of further coin tosses would not change the outcome.

The particular probability law (the so-called standard exponential distribution) used to generate the random number *Z* in this problem is not crucial. Any probability distribution over the positive real numbers that assigns positive probability to any interval of positive length does the job.

---

[87] Ok

[88] Yes, but without information on the potential amounts, how can one demonstrate that the probability of leaving the gain with the largest amount is $p = 1/2 + ε$ with $ε > 0$ and ε fixed in advance? In fact this approach is theoretically correct but actually unusable.

[89] How can one simulate a number whose probability of occurrence of any value is exactly equal to zero?

[90] How can one show that this probability is greater than 1/2 + *ε*, for *ε* > 0 fixed?

This problem can be considered from the point of view of game theory, where we make the game a two-person zero-sum game with outcomes win or lose, depending on whether the player ends up with the higher or lower amount of money. The organiser chooses the joint distribution of the amounts of money in both envelopes, and the player chooses the distribution of *Z*. The game does not have a "solution" (or saddle point) in the sense of game theory. This is an infinite game and von Neumann's minimax theorem does not apply.[31] [91]

**History of the paradox**

The envelope paradox dates back at least to 1953, when Belgian mathematician Maurice Kraitchik proposed a puzzle in his book Recreational Mathematics concerning two equally rich men who meet and compare their beautiful neckties, presents from their wives, wondering which tie actually cost more money. It is also mentioned in a 1953 book on elementary mathematics and mathematical puzzles by the mathematician John Edensor Littlewood, who credited it to the physicist Erwin Schroedinger. Martin Gardner popularized Kraitchik's puzzle in his 1982 book Aha! Gotcha, in the form of a wallet game:

Two people, equally rich, meet to compare the contents of their wallets. Each is ignorant of the contents of the two wallets. The game is as follows: whoever has the least money receives the contents of the wallet of the other (in the case where the amounts are equal, nothing happens). One of the two men can reason: "I have the amount *A* in my wallet. That's the maximum that I could lose. If I win (probability 0.5), the amount that I'll have in my possession at the end of the game will be more than 2*A*. Therefore the game is favourable to me." The other man can reason in exactly the same way. In fact, by symmetry, the game is fair. Where is the mistake in the reasoning of each man? [92]

In 1988 and 1989, Barry Nalebuff presented two different two-envelope problems, each with one envelope containing twice what's in the other, and each with computation of the expectation value 5*A*/4. The first paper just presents the two problems, the second paper discusses many solutions to both of them. The second of his two problems nowadays the most common, and is presented in this article. According to this version, the two envelopes are filled first, then one is chosen at random and called Envelope A. Martin Gardner independently mentioned this same version in his 1989 book Penrose Tiles to Trapdoor Ciphers and the Return of Dr Matrix. Barry Nalebuff's asymmetric variant, often known as the Ali Baba problem, has one envelope filled first, called Envelope *A*, and given to Ali. Then a fair coin is tossed to decide whether Envelope *B* should contain half or twice that amount, and only then given to Baba. [93]

---

[91] This approach is entirely theoretical and is useful only if the probability law of the amounts in the envelopes is at least partially known. It is a disguised Bayesian approach.

[92] Person *A* cannot, at the same time, estimate that he has one chance out of two of winning and estimate that the portfolio of person *B* contains more than his own. The same goes for person *B*. This is the error that led to saying that this game would be "win-win". It's another problem...

[93] The Ali Baba alternative is very interesting; the question is whether it is worthwhile for the two protagonists to swap their envelopes. If they do, one will end up with double the sum he owns and the other will end up with half of what he owns. It is clear that Ali has no reason to think he is more likely to end up with double than with half of his sum. This is also the case for Baba. Let *X* be the value discovered by Ali and *Y* value discovered by Baba. Ali believes rationally that its conditional expectation by swapping is E(Baba|Ali=*X*) = 5/4*\**X*. Baba does the same: E(Ali|Baba=*Y*) = 5/4*\**Y*. The reasoning of each is correct and thus leads them to swap their envelopes. To summarize: 1) Ali=*X* ; 2) Baba=*Y*, P(*Y*=2*X*)=P(*Y*=*X*/2)=1/2, P(*X*=*Y*/2)=P(*X*=2*Y*)=1/2 ; 3) Ali+Baba=*X*+*Y* ; 4) For Ali, E(Baba|Ali=*X*) = 5/4*\**X* ; 5) For Baba, E(Ali|Baba=*Y*)=5/4*\**Y* ; 6) E(Baba|Ali=*X*) + E(Ali|Baba=*Y*) = 5/4*\*(*X*+*Y*).

What seems shocking is point (6). However, this calculation does not make sense in this problem. It adds two conditional expectation estimates that are not derived from the same information. Neither of the two protagonists possesses the same information. Ali knows that his envelope

## Notes and references

---

contains *X* and Baba knows that his envelope contains *Y*. There is nothing inconsistent or irrational that each man should produce a conditional estimate expectation whose sum is 5/4 *(X + Y).

*******************************************************************

**Conclusion**

The two envelopes paradox is remarkable. It highlights the need to differentiate the concepts of "whatever *x* known" and "whatever *x* unknown ". Indeed, it is expected that these two concepts are equivalent but the corresponding models differ fundamentally. While it is assumed that *x* is known (whatever *x*), *x* is a variable that does not need to be probabilized. However, when *x* is unknown, we have to provide a probability distribution for *x*. Yet in the latter case, there is no uninformative probability law for the amounts in the envelopes.

Ali's alternative is very instructive. In this version, the gain expectancy estimates are individually rational and there is no need or reason for the sum of these estimates to be equal to the sum of the two envelopes.

In the original two envelopes paradox, the two cases are as follows:

- If the first envelope is opened and the amount discovered is *x*, the estimate of the expected gain in swapping is possible, coherent, rational, and its value is 5/4*x*. Then, you should swap.

- If the first envelope is not opened, the only probabilistic reasoning that does not use supplementary information leads to estimating expectations as infinite amounts of each envelope. This means therefore that there is no point in swapping envelopes.

To give a frequentist interpretation of estimating the conditional expectation by swapping at 5/4*x*, we must understand what type of repeated experience it corresponds to. It is the schema where the animator chooses amount x for an envelope, then selects with probability 1/2 the amount *x*/2 or 2*x* the amount for the other envelope. Then, you randomly select an envelope and you only study all realizations where the first selected envelope contains amount x. Also, conditional on the fact that the first envelope contains value *x*, we gain more, on average, by swapping.

When everybody is trying to demonstrate unilaterally that the only rational and coherent analysis is irrational and incoherent, we obtain a paradox…